\documentclass [a4paper,russian]{article}
\usepackage[cp1251]{inputenc}
\usepackage{amsmath, amsfonts, amssymb}
\usepackage[russian,english]{babel}
\usepackage[T2A]{fontenc}
\usepackage{babel,indentfirst,epigraph}
\usepackage{graphicx,graphics,verbatim}

\begin{document}
\vspace*{0.15in}
\centerline{
\large
\textbf{THE TRANSMUTATION OPERATORS}
}
\centerline{
\large
\newline\textbf{AND CORRESPONDING HYPERBOLIC EQUATIONS.}
}
\textbf\newline{}
\centerline{
V.I. MAKOVETSKY\textsuperscript{*},
}
\centerline{
* Makov-Vikror-Sakh@yandex.ru}

\section*{Abstract}
For arbitrary second-order differential operators, the existence conditions and the construction of intertwining transmutation operators are shown. In an explicit form found hyperbolic equations with two independent variables and their solutions leading to the kernels of the transmutation operators.

\section{Definitions and preliminary comments.}
By analogy with Volterra integral equations, we introduce the \emph {operators transformations of the first kind} ($OP_{I}$)
\begin{equation}\label{1.1}
Tf(x) = f_1(x) = \int\limits_0^x K(x,t)f_0(t)dt
\tag{1.1} \end{equation}
и \emph{operators transformations of the second kind} ($OP_{II}$)
\begin{equation}\label{1.2}
Tf(x) = f_1(x) = f_0(x)+\int\limits_0^x K(x,t)f_0(t)dt
\tag{1.2} \end{equation}
The transfer function $ K (x, t) $ is called the \emph {kernel} of the transformation operator.

If $ K (x, x) = \gamma \neq 0 $ in (\ref {1.1}), then by differentiation the equation (\ref {1.1}) goes to (\ref {1.2})
\begin{equation}\label{1.3}
 \Psi(x)=f'_1(x) = \gamma f(x) + \int\limits_0^x K(x,t)f(t)dt
\tag{1.3} \end{equation}
Due to this fact, only ($OP_{I}$) will be studied later.

\textbf{Comment.} The singularities of the kernel and the coefficients of the subsequent differential equations involved in the construction of $K (x, t)$ require more correct recording of the proposed definition. Exactly
\begin{equation}\label{1.4}
Tf(x) = f_1(x) = \int\limits_\varepsilon^{x-\delta} K(x,t)f(t)dt
\tag{1.4} \end{equation}
with the subsequent passage $\varepsilon \rightarrow 0$ and $\delta \rightarrow 0 $. These refinements will be correctly spelled out when setting the conditions imposed on the kernel of the conversion operator.

\section{The spectral problem and the basic identity for the transmutation operator.}
The most interesting part of the transformation operator is a specific kind of its kernel. There are many approaches to its construction, but we will choose a variant related to the spectral properties of the accompanying operators. Namely, if (\ref {1.4}) is satisfied, eigenvalue problems for the basic and transformed functions will be considered.
\begin{align}\label{2.1}
  & A f_0 = \lambda f_0 \tag{2.1a} \\
  & B f_1 = \lambda f_1 \tag{2.1b}
\end{align}
Changes in the right side of the second line lead to the result
\begin{equation*}
B f_1 = \lambda f_1 = \lambda(T f_0) = T(\lambda f_0) = T(A f_0)
\end{equation*}
In the left part
\begin{equation*}
  B f_1 = B(T f_0)
\end{equation*}
So
\begin{equation*}
B(Tf_0) = T(Af_0)
\end{equation*}
Due to the arbitrariness of the function $f_0(x)$, the \emph {main identity for the transmu\-tation operator} come into being
\begin{equation}\label{2.2}
BT = TA
\tag{2.2}\end{equation}
Note that the derivation of the identity did not establish a specific type of operators A and B. It is important only the existence of their eigenvalues and eigenfunctions.

\section{Formulation of the problem leading to hyperbolic equations.}

Let A and B be differential operators of the second order, defined by differential expressions and boundary conditions on the half-line $(0, + \infty)$
\begin{equation}\label{3.1}
\left \{
\begin{aligned}
    & B  = a_1(x)\partial_{x,x}(\circ) +b_1(x)\partial_x (\circ) +c_1(x) (\circ); \\
    & \partial \left.(\circ)\right|_{x=0} -h_1\left.(\circ)\right|_{x=0} =0    \\
    \end{aligned}
    \right.
\tag{3.1}\end{equation}
and
\begin{equation}\label{3.2}
\left \{
\begin{aligned}
    & A  = \partial_t(a_0(t)\partial_t (\circ)) +\partial_t (b_0(t)(\circ)) +c_0(t) (\circ); \\
    & \partial \left.(\circ)\right|_{t=0} -h_0\left.(\circ)\right|_{t=0} =0     \\
\end{aligned}
    \right.
\tag{3.2}\end{equation}
The problem of constructing the kernel of the transformation operator K(x, t) is solved only for the variant
\begin{equation*}
  a_1(x)= a_0(t) = -1; \quad b_1(x) = b_0(t) = 0;
\end{equation*}
However, there are very few spectral problems with such a classical form (\cite{Everitt2004},\cite{Titchmarsh1960}). Therefore, we investigate the construction of K (x, y) for arbitrary coefficients in the operators A and B. It turns out that the kernel K (x, t) is a solution of linear hyperbolic equations with two independent variables (\cite {Courant1964}, Chapter V). This fact allows a) to construct the transformation operators if the integrals of hyperbolic equations are known (in particular, the Riemann function), b) to find solutions of hyperbolic equations if the transformation operators are obtained by other means, for example, in the study of integral tables.


\subsection{Overview of methods for constructing conversion $\\$ operators for the classical equations}

In general, there are two approaches to generating the kernel K (x, t) for a pair of operators
\begin{equation}\label{3.1.1}
B = -\partial_{xx}(\circ) + c_1(x)(\circ); \quad
A = -\partial_{tt}(\circ) + c_0(t)(\circ);
\tag{3.1.1} \end{equation}
under the corresponding initial conditions on the half-line $(0,+\infty)$

The first direction is a direct substitution of transmutation (\ref {1.1}) into the basic identity (\ref {2.2}). It is developed in the works (\cite {GelfandLevitan1951}, \cite {Levitan1984}, \cite {Naimark1969}, \cite {LionsJL1956}) for inverse Sturm-Liouville problems. The most interesting is that, despite the direct need to include transformation operators in this theory, the specific form of the kernel K (t, x) of the creators of the concept was of an extremely low interest. In the book (\cite{Marchenko1977}, p. 19) 'The transformation operators \dots and their inverses play a very important role in the Sturm-Liouville spectral theory. To solve many of the basic problems of this theory, one fact of their existence is sufficient.".\- Great importance was attached to their estimates and asymptotic expansions, but not to the immediate construction of the kernel K (x, t). Conversely, when solving the inverse Sturm-Liouville problem, the core itself satisfies the Gel'fand-Levitan integral equation
\begin{equation}\label{3.1.2}
\chi(x,t) + \int\limits_0^x K(x,s)\chi(s,t)ds +K(x,t) = 0
\tag{3.1.2} \end{equation}
with initial function
\begin{equation}\label{3.1.3}
\chi(x,y) = \frac{\partial^2}{\partial x \partial y}
\int\limits_{-\infty}^{+\infty}\frac{\sin \sqrt{\lambda x}\sin \sqrt{\lambda y} }{\lambda} d \tau(\lambda)
\tag{3.1.3} \end{equation}
and for $c_1 (x) = 0; \: c_2 (x) = q (x)$ the coefficient q (t) is found by K (x,t) using the relation
\begin{equation}\label{3.1.4}
q(x) = 2 \partial_x K(x,x)
\tag{3.1.4} \end{equation}
This is the famous problem of determining the coefficient of the differential operator $ A = - \partial_{xx} (\circ) + q (x) (\circ) $ by means of the spectral function $ \tau (\lambda) $, on the assumption of existence $Af = \lambda f$ (\cite{GelfandLevitan1951}).

The second direction is to obtain the kernel K (x, t) by solving a hyperbolic equation, especially using the Riemann function. For the Sturm-Liouville problem, this is the Marchenko method(\cite{Marchenko1977}, ch. I). For generalized translation operators, this approach developed in the works(\cite{Levitan1949}, \cite{CarrolShowalter1976}, \cite{Carrol1985}, \cite{Trimeche1988}). It naturally fits into the definition of the shift operator $ T^s_t $, given in Chapter I of the book B. \, M. Levitan (\cite{Levitan1973}). Namely, for a space of variables $\Omega$ and a function space $ C = C (\Omega) $ each function $ f (t) \in C $ is associated with a function of two points $(t, s)$ such that
\begin{equation}\label{3.1.5}
F(s,t) = T^s_t f(t)
\tag{3.1.5} \end{equation}
But such an operation is just performed by the Riemann functional, which the initial conditions
\begin{equation*}
u(x,0)=\phi(0); \quad  \partial_t \phi(x,0) = \psi(x)
\end{equation*}
converts to the value of u (x, t) at some point(\cite{Marchenko1977}, Ch. 1, \S\,1)
\begin{equation}\label{3.1.6}
u(x_0,y_0) = T^{t_0}_{x_0}[\phi(x),\psi(x)]
\tag{3.1.6} \end{equation}
or get rid of one of the functions for complete coincidence with the shift operator, it is possible in one of two ways: either by zeroing out any of them, or (as in \cite {Marchenko1977}, Ch. 1)
\begin{equation*}
  \psi(x) = \partial_x \phi(x)
\end{equation*}
Particularly successful, Riemann's immediate approach turns out to be in various modifi\-cations of the generalized hyperbolic Euler-Poisson-Darboux equation (\cite {Levitan1951}). (\cite {Sitnik2016}, ch. II, Kopson's lemma),(\cite{Sitnik2017}), (\cite{ShishkinaSitnik2017}).

Note that the transformation (\ref {3.1.6}) is not a privilege of only hyperbolic equations. In spite of the fact that the interweaving Vekua-Erdei-Lowndens operator carrying out a shift in the spectral parameter was developed in a hyperbolic variant, its original formulation was based on the Laplace equation (that is, referred to the elliptic type of equations)(\cite{Vekua1948}, Ch. 1)

Beyond our article, the fruitful direction of compositional construction of transforma\-tion operators, when several transmutations are laid in the basis of intertwining pairs, and further types of OP-operators are obtained by applying them consistently to the original equations (\cite{RyzhkovaSitnik2016}, \cite{ShishkinaSitnik2017A}, \cite{FitouhiShishkina2018}). This is due to the fact that we are interested in the methods of constructing the initial, seed nuclei of transmutation, and not the iterative methods of their influence and recalculation.

\section{Conditions for constructing the kernel of the transmutation operator.}
Let the transformation operator be realized in the form (\ref {1.4}) under the condition $ \varepsilon \rightarrow 0 $ and $ \delta \rightarrow 0 $. In order for it to transform the solution of the differential operator (\ref {3.1}) into the solution of the differential operator (\ref {3.2}) on the semiaxis $ (0, + \infty) $, it is necessary to satisfy the following conditions: \\
\newline
\emph{a) The kernel of the operator satisfies the hyperbolic equation}
\begin{align}\label{4.a}
  & \partial_t\{a_0(t)\partial_tK(x,t)\}-b_0(t)\partial_t K(x,t) +c_0(t)K(x,t)= \notag\\
  & = a_1(x)\partial_{xx}K(t,x)+b_1(x)\partial_x K(x,t) +c_1(x) K(x,t) \notag
\tag{4.a}\end{align}
 \emph{b) On the assumption of $t\rightarrow x-\delta$ } и  $\delta \rightarrow 0$
\begin{align}\label{4.b}
  & \emph{b1)}\quad a_0(x) = a_1(x) = a(x) \notag\\
  & \emph{b2)}\quad 2 a(x)\frac{d K(x,x-\delta)}{dx} + (b_1(x)-b_0(x))K(x,x-\delta)=0
\tag{4.b}\end{align}
\emph{c) On the assumption of $t = \varepsilon\rightarrow 0 $}
\begin{equation}\label{4.c}
\left\{a(\varepsilon)\,\left[\partial_t K(x,t)\right]_{t=\varepsilon}
- b_0(\varepsilon)K(x,\varepsilon)
-h_0 a(\varepsilon)K(x,\varepsilon)\right\}f_0(\varepsilon) \rightarrow 0
\tag{4.c}\end{equation}
\emph{d) Condition at the vertex} On the assumption of $ \delta < \varepsilon, \: \delta \rightarrow 0, \: \varepsilon \rightarrow 0$
\begin{equation}\label{4.d}
K(\varepsilon,\varepsilon-\delta)*f_0(\varepsilon) \rightarrow 0;
\tag{4.d}\end{equation}
\newline
Conversely, the fulfillment of all the above conditions is sufficient for the existence of a transformation operator of the form (\ref {1.4}). Note that in points c) and d) it is not necessary to know the explicit form of the function $ f_0 (x) $. Only the rate of $ f_0 (\varepsilon) $ aspiration is important to zero when $ \varepsilon \rightarrow 0 $ to compensate for the singularities of the coefficients $ a, b, c $ at zero.

\subsection{Theorem proving}
We prove the assertion of the theorem, generalizing the method(\cite{Levitan1984},\cite{Naimark1969}). For ease, we introduce the notation
\begin{equation}\label{4.1.1}
  K_0(x) = K(x,x-\delta); \quad f(x) = f_0(x);
\tag{4.1.1}\end{equation}
First, the operation $TA$ will be performed.
\begin{align}\label{4.1.2}
  & TA(f(x))= \rho\{\partial_x(a_0(x)\partial_x f(x))+
  \partial_x(b_0(x)f(x)+c_0(x)f(x)\}+ \notag\\
  & +\int_{\varepsilon}^{x-\delta}K(x,t)\{\partial_x(a_0(x)\partial_x f(x))+
  \partial_x(b_0(x)f(x))+c_0(x)f(x)\}dt\notag
\tag{4.1.2}\end{align}
The integral with the first term is taken twice in parts. Similarly, by parts, only the second term will be transformed once. This leads to the following result
\begin{align*}
  & \text{TA}(f(x))= \rho\{\partial_x(a_0(x)\partial_x f(x))+
  \partial_x(b_0(x)f(x)+c_0(x)f(x)\}+a_0(x)K_0(x)\partial_x f(x)+\\
  & + \left[K_0(x)b_0(x)- a_0(x)\left.\{\partial_t K(x,t)\}\right|_{t=x-\delta}\right]f(x)
  - a_0(\varepsilon)K(x,\varepsilon)\left.\{\partial_t f(t)\}\right|_{t=\varepsilon}+\\
  &+a_0(\varepsilon)\left.\{\partial_t K(x,t)\}\right|_{t=\varepsilon}f(\varepsilon)
  -b_0(\varepsilon)K(x,\varepsilon)f(\varepsilon)+\\
  & +\int_{\varepsilon}^{x-\delta}\left\{ \partial_t [a_0(t)\partial_t K(x,t) ]
  -b_0(t)\partial_t K(x.t) + c_0(t)K(x,t)\right\} f(t)dt
\end{align*}
The second on the order is the operation $BT$
\begin{align}\label{4.1.3}
B(Tf(x)) = \left\{
a_1(x)\partial_{x,x}(\circ) +b_1(x)\partial_x (\circ) +c_1(x) (\circ)
\right\}
\left\{
\rho f(x) + \int\limits_\varepsilon^{x-\delta} K(x,t)f(t)dt
\right\}
\tag{4.1.3}\end{align}
Calculating the integral over the variable upper limit leads to the result
\begin{align*}
  & B(Tf(x)) = \rho\{
 a_1(x)\partial_{x,x}f(x) +b_1(x)\partial_x f(x) +c_1(x) f(x)
  \}+ \\
  & +\{a_1(x)\partial_x K_0(x)+a_1(x)\left.[\partial_x K(x,t)\rfloor\right|_{t=x-\delta}+
  b_1(x) K_0(x)\}f(x)+a_1(x)K_0(x)\partial_x f(x)+\\
  & + \int_{\varepsilon}^{x-\delta}K(x,t)
  \{a_1(x)\partial_{xx}K(x,t)+b_1(x)\partial_x K(x,t) + c_1(x)K(x,t)\}f(x)
\end{align*}
A comparison of the integrands leads to (\ref {4.a}). In view of the arbitrariness of $f(x)$, the coefficients of the function and its first derivative must vanish separately. Comparison of the elements before the first derivative yields (\ref {4.b}.1). If we take this fact into account in the coefficient adjacent to $f(x)$, and use equality
\begin{equation*}
\frac{d K(x,x-\delta)}{dx} = \left.[\partial_x K(x,t)\rfloor\right|_{t=x}
+ \left.[\partial_t K(x,t)\rfloor\right|_{t=x}
\end{equation*}
then the grouping of elements before $f(x)$ establishes a match (\ref {4.b}.2). It remains to group the string of initial conditions when $ t = \varepsilon \rightarrow 0 $. All its elements are entirely in the $ TA $ operator. The expression to take place
\begin{equation*}
-a_2(\varepsilon)K(x,\varepsilon)\left.\{\partial_t f(t)\}\right|_{t=\varepsilon}+
a(\varepsilon)\left.\{\partial_t K(x,t)\}\right|_{t=\varepsilon}f(\varepsilon)
  -b_2(\varepsilon)K(x,\varepsilon)f(\varepsilon)
\end{equation*}

The result is fixed in the condition (\ref {4.c}). To formulate the condition at the vertex, we take the derivative of (\ref {1.4}).
\begin{equation}\label{4.1.4}
\partial_x Tf(x) = K(x,x-\delta)f(x-\delta) +
\int\limits_{\varepsilon}^{x-\delta}\partial_x K(x,t) f(t) dt
\tag{4.1.4} \end{equation}
In point $x=\delta+\varepsilon$
\begin{equation*}
\partial \left.Tf(x)\right|_{x=\delta+\varepsilon} = \rho*
\partial \left.f(x)\right|_{x=\delta+\varepsilon} +
K(\delta+\varepsilon,\varepsilon)*f(\varepsilon)
\end{equation*}
We take the initial conditions
\begin{equation*}
\partial \left.f_1(x)\right|_{x=\varepsilon} -h_1\left.f_1(x)\right|_{x=\varepsilon} = 0; \quad
\partial \left.f_0(x)\right|_{x=\varepsilon} -h_0\left.f_0(x)\right|_{x=\varepsilon} = 0;
\end{equation*}
The result is(\ref{4.d})

\subsection{An example of a general transmutation operator and the associative hyperbolic equation}

 The transmutation operator presented below(\cite{PrudnikovBrychkovMarichev2002}, vol II, № 2.21.2 (2))
\begin{align}\label{4.2.1}
& \qquad \qquad \qquad \qquad \qquad f_1(x) = \int\limits_0^x K(x,t) f_0(x)dx  \tag{4.2.1 a} \\
& f_1(x) = \frac{1}{2}B(\lambda-\beta,\beta)x^{2(\lambda-1)+p}C_p^{\lambda-\beta}(x); \quad
f_0(t) = t^{2(\lambda-\beta)+p-1}C_p^\lambda(t); \tag{4.2.1 b}\\
& \qquad \qquad \qquad \qquad \qquad K(x,t) = \left(x^2-t^2\right)^{\beta-1}; \tag{4.2.1 c}
\end{align}
by $p=2n+1, \; \varepsilon = 0 \;\emph{или}\;1, \; Re \beta > 0, \; Re(\lambda-\beta) > -n-\frac{\varepsilon}{2}$; $B(\lambda-\beta,\beta)$ - beta function, connects two Gegenbauer polynomials $f(x)=C_m^\alpha(x)$, satisfying the ordinary differential equation
\begin{equation*}
(1-x^2)\partial_{xx}f(x)-(2\alpha+1)x\partial_x f(x)+m(m+2\alpha)f=0
\end{equation*}
on the interval $[-1,1]$.

Let us find the differential equation for the product of the ultraspherical polynomial to an arbitrary power of the independent variable-in other words, the equation for
\begin{equation*}
\chi(x) = x^\delta C_m^\varrho(x)
\end{equation*}
Of course, $ \chi (x)/x^\delta$ is a solution of the original Gegenbauer equation, from which it follows that
\begin{align*}
  & x^2(x^2-1)\partial_{xx} \chi(x) + x(2\delta -x^2(2(\delta-\varrho)-1)\partial_x \chi(x) - \\
  & \qquad - [(m+\delta)(m-\delta+2\varrho)x^2+\delta(\delta+1)]\chi(x) = 0;
\end{align*}
By means of the presented equality, the following results can be established.
\newline
By $\delta = 2(\lambda-1)+p; \quad m = p; \quad \varrho = \lambda-\beta$
\begin{align*}
  & a_1(x) = x^2(x^2-1); \\
  & b_1(x) = x*[(5-2(p+\beta+\lambda))x^2+2(p+2(\lambda-1))]; \\
  & c_1(x) = 4(\beta-1)(p+\lambda-1)x^2-(p+2\lambda-2)(p+2\lambda-1)
\end{align*}
the ordinary differential equation $ B f_1 (x) = 0 $ is satisfied (see \ref {3.1}). The last equality is equivalent to the search for an eigenfunction for the zero-characteristic number.

If $\delta = 2(\lambda-\beta)+p-1; \quad m = p; \quad \varrho = \lambda$
\begin{align*}
  & a_0(t) = t^2(t^2-1); \\
  & b_0(t) = -t*[(2p-4\beta+2\lambda+1)t^2+2(2(\beta-\lambda)-p)]; \\
  & c_0(t) = 4(\beta-1)(\beta-\lambda-p-1)t^2 + [4\beta p +2\beta -p(1+p)-4(\beta-\lambda)^2 -2\lambda-4\lambda p]
\end{align*}
then execute $A f_0(t) = 0$ (\ref{3.2}).

The presented coefficients and the kernel K (x, t) completely satisfy (\ref{4.a})

\section{Transmutations operators with kernels Erdelyi - Kober type.}

We consider the Erdelyi - Kober kernels in the original formulation, as well as their generalization.
\subsection{Traditional fractional Erdelyi - Kober derivatives}
The usual Erdelyi - Kober operator is defined as(\cite{SamkoKilbasMarichev1987}, Ch. 4))
\begin{equation}\label{5.1.1}
I^{\alpha}_{a+,\sigma,\eta}f(x) = \frac{\sigma x^{-\sigma(\alpha+\eta)}}{\Gamma(\alpha)}
\int\limits_a^x (x^\sigma - t^\sigma)^{\alpha-1}t^{\sigma\eta+\sigma-1}f(t)dt
\tag{5.1.1}\end{equation}
where $ \Gamma (\alpha) $ is the natural gamma function. One of the variants of the Erdelyi - Kober operator is presented in Section 4.2.

Despite of the apparently simple form of the kernel
\begin{equation}\label{5.1.2}
K(x,t) = \frac{\sigma x^{-\sigma(\alpha+\eta)}}{\Gamma(\alpha)}
(x^\sigma - t^\sigma)^{\alpha-1}t^{\sigma\eta+\sigma-1}
\tag{5.1.2}\end{equation}
finding the explicit form of a hyperbolic equation that it can satisfy is an extremely complicated problem. Therefore, we consider the narrowing of a group of parameters.

\emph{Example 5.1.1. The connection between the Poisson transformation for the Bessel function with the hyperbolic equation.}
\begin{equation}\label{5.1.3}
K(x,t) = \frac{1}{2^{\nu-1}*\Gamma\left(\nu+\frac{1}{2}\right)}*
\frac{(x^2-t^2)^{\nu-\frac{1}{2}}}{x^{{\nu-\frac{1}{2}}}}
\tag{5.1.3}\end{equation}
The kernel is the solution of the hyperbolic equation
\begin{equation}\label{5.1.4}
\partial_{tt}K(x,t)+K(x,t)=
\partial_{xx}K(x,t)+\left[
1-\frac{\nu^2-\left(\frac{1}{2}\right)^2}{x^2}
\right]K(x,t)
\tag{5.1.4}\end{equation}
Recording the coefficients in (\ref{4.a})
\begin{equation*}
a(x) = 1; \quad b_1(x) = b_0(t) = 0; \quad
c_1(x) = 1-\frac{\nu^2-\left(\frac{1}{2}\right)^2}{x^2};\quad c_0(t) = 1;
\quad \rho = 0;
\end{equation*}
Of course, the unit in the free terms could be omitted, but it turns out to be useful in the study of the spectral problem (\ref {2.1}). It is easy to see that the functions
\begin{equation*}
\quad f_1(x) =  \sqrt{x}J_\nu (x); \quad f_0(x) = \cos(x);
\end{equation*}
are solutions of (2.1), with operators(\ref{3.1},\ref{3.2}) with eigenvalue $\lambda = 1$. The initial condition for $\varphi_0(x)$ - is $h_0 =0$. As for $f_1(x)$, the asymptotics
\begin{equation*}
\partial \varphi_1(x) - h_1 \varphi_1(x)
\end{equation*}
in the neighborhood of zero becomes
\begin{equation*}
\frac{1}{2^\nu \Gamma(\nu+1)}
\left[x^{\nu-\frac{1}{2}}\left(\nu+\frac{1}{2}\right)
+h_1 x^{\nu+\frac{1}{2}}
\right]
\end{equation*}

It is clear that when $\nu > \frac{1}{2}$ the value of the initial constant can be any, including $h_1 = 0$. If $x \neq 0$, then $\delta \rightarrow 0$ The remaining conditions of the items (\ref {4.a} - \ref {4.d}) are correct.

The equality
\begin{equation}\label{5.1.5}
\sqrt{x}J_\nu (x) =  \frac{1}{2^{\nu-1}\sqrt{\pi}\Gamma\left(\nu+\frac{1}{2}\right)}
\int_0^x \frac{(x^2-t^2)^{\nu-\frac{1}{2}}}{x^{{\nu-\frac{1}{2}}}}\,\cos(t)dt
\tag{5.1.5}\end{equation}
defines the Erdelyi - Kober fractional derivative (\ cite {Ross1975}) and, simultaneously, is a Poisson representation for the Bessel function. To verify the latter, it suffices to make two substitutions $ t = x * \xi $ and $ \xi = \sin (\theta) $. The resulting formula (\cite{BatemanErdelyi1974})
\begin{equation}\label{5.1.6}
\Gamma\left(\nu+\frac{1}{2}\right)J_\nu(x) =
\frac{2}{\sqrt{\pi}}\left(\frac{x}{2}\right)^\nu
\int\limits_0^{\frac{\pi}{2}}\cos(x\sin\theta)(\cos\theta)^{2\nu}d\theta
\tag{5.1.6}\end{equation}
was obtained by S.D. Poisson in 1823 for $ \nu> - \frac {1}{2} $. In our variant, the restriction $\nu> + \frac {1}{2} $ is severer, but this allowed us to write the hyperbolic equation (\ref {5.1.4}) for the kernel of the transformation operator(\ref{5.1.3})(\cite{Makovetsky2017})

\emph{Example 5.1.2. The Erdelyi - Kober kernel, which leads to the first Sonin integral}.

For real parts $ Re(\mu)> -1 $ and $ Re(\nu)> -1 $, the first Sonin integral (\cite {Watson1949}, Ch. 12) connects Bessel functions with different indices.
\begin{equation}\label{5.1.7}
J_{\mu+\nu+1}(x) = \frac{x^{\nu+1}}{2^\nu \Gamma(\nu+1)}
\int\limits_0^{\frac{\pi}{2}}J_\mu(x\sin(\theta)\sin^{\mu+1}\theta
\cos^{2\nu+1}\theta d\theta
\tag{5.1.7}\end{equation}
The replacements $ \xi = \cos t $ and $ t = x \setminus \xi $ and the rearrangement of the elements are given by the integral transmutation operator of the first kind (\ ref {1.1})
\begin{align}\label{5.1.8}
& f_1(x) = x^{\mu+1}J_{\mu+\nu+1}(x); \quad
f_0(t) = t^{\mu+1}J_\mu(t); \tag{5.1.8 a}\\
& \qquad  K(x,t) =
\frac{1}{2^\nu \Gamma(\nu+1)}\frac{ \left(x^2-t^2\right)^\nu}{x^\nu}; \tag{5.1.8 b}
\end{align}
Coefficients
\begin{align*}
  & a_1(x)= 1; \quad a_0(t) = 1; \\
  & b_1(x) = -\frac{2\mu+1}{x}; \quad b_0(t) = -\frac{2\mu+1}{t}; \\
  & c_1(x) = 1-\frac{\nu(\nu+2\mu+2)}{x}; \quad c_0(t) = 1;
\end{align*}
guarantee the execution of $A f_0 (t) = 0; B f_1 (x) = 0 $. If each of the coefficients of the free term is reduced by one, then similar operators lead to an eigenvalue $ \lambda = 1 $. Conditions (\ ref {4.a} - \ ref {4.d}) are also fully realized. Here, the equality $ b_1 (x) = b_0 (x) $ is especially helpful. The similarity of the kernels of Examples 1 and 2 indicates the general properties of the Sonin and Poisson integrals.


\emph{Example No. 5.1.3 of the Erdelyi - Kober core with a smooth kernel.}

Previous images of nuclei have a characteristic feature at the point $ x = 0 $. All this is connected with the nature of the Bessel equations, which even for spectral decomposition requires the allocation of subdomains by a real axes(\cite{Naimark1969} \textsection \, 21). The kernel of the form
\begin{equation}\label{5.1.9}
K(x,t) = t(x^2-t^2)^{\beta-1}
\tag{5.1.9}\end{equation}
transforms a sine into a Bessel function of an arbitrary half-integer order (\cite{PrudnikovBrychkovMarichev2002}, vol I, 2.5.55 (7))
\begin{equation}\label{5.1.10}
\int\limits_0^x K(x,t) \sin(\omega t) dt = \frac{\sqrt{\pi}}{2}x
\left(\frac{2x}{\omega}\right)^{\beta-\frac{1}{2}}\Gamma(\beta
)J_{\beta+\frac{1}{2}}(\omega*x)
\tag{5.1.10}\end{equation}
Direct calculation confirms the equality
\begin{equation}\label{5.1.11}
\partial_{tt}K(x,t)+\omega^2 K(x,t)=
\partial_{xx}K(x,t)-\frac{2\beta}{x}\partial_x K(x,t) + \omega^2 K(x,t)
\tag{5.1.11}\end{equation}
Its coefficients
\begin{equation*}
a(x) = 1; \quad b_1(x) = -\frac{2\beta}{x}, \quad b_0(t) = 0; \quad
c_1(x) = \omega^2;\quad c_0(t) = \omega^2;
\end{equation*}
Right part
\begin{equation*}
f_1(x) = \frac{\sqrt{\pi}}{2}x
\left(\frac{2x}{\omega}\right)^{\beta-\frac{1}{2}}\Gamma(\beta
)J_{\beta+\frac{1}{2}}(\omega*x)
\end{equation*}
satisfies the spectral equation for an ordinary second-order differential equation
\begin{equation*}
\partial_{xx}f_1(x)-\frac{2\beta}{x}\,\partial_x f_1(x)+\omega^2 f_1(x)= 0
\end{equation*}
for the conventional value $\lambda = -\omega^2$

\emph{Example 5.1.4. Integration into classes of integrals with equivalent transmutation kernels}

The similarity of the transformation kernels makes it possible to distinguish classes of integrals with common properties. For an example, let us return to the Poisson integral (\ref {5.1.5}) by replacing it
\begin{equation*}
\nu \rightarrow \beta - \frac{1}{2}; \quad x \rightarrow \omega x; \quad t \rightarrow \omega t;   \end{equation*}
The result will be(\cite{PrudnikovBrychkovMarichev2002}, vol I, 2.5.6 (1)) if  $ x > 0, \, Re(\beta) > 0; \, |arg \omega | < \pi$
\begin{equation*}
\frac{\sqrt{\pi}}{2}  \left(\frac{2x}{\omega}\right)^{\beta-\frac{1}{2}}
\Gamma(\beta)J_{\beta-\frac{1}{2}}(\omega x)=
\int\limits_0^x (x^2-t^2)^{\beta-1}\cos(\omega t)dt
\end{equation*}
Proceeding from the identity for the modified Bessel function
\begin{equation*}
  I_\alpha(x) = i^{-\alpha}J_\alpha (ix)
\end{equation*}
an equivalent integral appears (\cite{PrudnikovBrychkovMarichev2002}, vol I, 2.4.3 (7))
\begin{equation*}
\frac{\sqrt{\pi}}{2}  \left(\frac{2x}{\omega}\right)^{\beta-\frac{1}{2}}
\Gamma(\beta)I_{\beta-\frac{1}{2}}(\omega x)=
\int\limits_0^x (x^2-t^2)^{\beta-1}\cosh(\omega t)dt
\end{equation*}

\subsection{Generalized fractional Erdelyi - Kober derivatives}
The power function built into the construction of Erdei-Kober's ordinary fractional derivative does not cover the full range of problems considered in applications. Hence, it is advantageous to replace it with more complex compo\-nents (\cite {Sneddon1975}). The need for operators of a new type proved to be so necessary that it began to develop immediately in several directions, beginning with equations of a mixed type of higher orders, up to singular wave equations and problems of a multidimensional random walk. In the list of literature, we presented only a small part of the work on these issues (\cite{Sprinkhuizen1979} - \cite{Karimov2018}). Historical reviews are displayed in articles (\cite{Colton1979},\cite{Karimov2017}). Very often the starting point for the derivation of fractional Erdei-Kober derivatives is the generalized Axially Symmetric Potential Equation (GASPE) or Generalized Axially Symmetric Potential Theory (GASPT), and the final product is associated with the application of the modified Euler-Poisson-Darboux equation. The latter is derived in each concrete case, and here it arises from the spectral approach to the theory of transformation operators (\ref {4.a} - \ref {4.d}).

\emph{Example 5.2} Fractional derivatives of the form
\begin{equation}\label{5.2.1}
If(x) = \int\limits_0^x K(x,t)f(t)dt
\tag{5.2.1}\end{equation}
where the kernels is
\begin{equation}\label{5.2.2}
K(x,t) = \sinh(\mu\sqrt{x^2-t^2)}\quad \emph{или} \quad
K(x,t) = \sin(\mu\sqrt{x^2-t^2})
\tag{5.2.2}\end{equation}
Obviously for small values of $\mu$ they will be converted into traditional fractional derivative. The general theory for them is not presented in the literature. Examine some equivalent groups. We are considering the transformation operator(\cite{PrudnikovBrychkovMarichev2002}, vol I, 2.4.8 (1)) in which
\begin{align}\label{5.2.3}
 & K(x,t) = \sinh(\mu\sqrt{x^2-t^2)} \tag{5.2.3 a}\\
 & f_1(x) = \frac{\pi}{2}\frac{\mu x}{\sqrt{\beta^2+\mu^2}}I_1(\sqrt{\beta^2+\mu^2} x);\quad f_0(t) = \cosh(\beta t); \tag{5.2.3 b}
\end{align}
With coefficients
\begin{align*}
  & a_1(x)= 1; \quad a_0(t) = 1; \\
  & b_1(x) = -\frac{1}{x}; \quad b_0(t) = 0; \\
  & c_1(x) = -(\beta^2+\mu^2); \quad c_0(t) = -\beta^2;
\end{align*}
we have the equalities $ Bf_1 (x) = 0 $ and $ Af_0 (t) = 0 $, which is equivalent to the eigenvalue problems $ \lambda = \beta^2 $. Hyperbolic equation and other conditions
(\ ref {4.a} - \ ref {4.d}) are satisfied because $ K (x, x) = 0 $. Replacements of the type $\beta \rightarrow \imath \beta $ and $ \mu \rightarrow \imath \mu $ in any combination lead to a whole group of equivalent by properties and to the kernel transformation operators (\cite{PrudnikovBrychkovMarichev2002}, vol I, 2.5.25; 2.5.53). We note that the presented operators can be derived from generalized derivatives related to Bessel functions as a kernel (they are considered below), but the simplicity of their construction deserves special attention.


\emph{Example 5.3 Lowndes fractional derivatives}
\begin{equation}\label{5.3.1}
I_{\mu,\gamma,\beta}f(x) = \int\limits_0^x
(x^2-t^2)^{\frac{\mu}{2}}J_\mu(\beta\sqrt{x^2-t^2})t^{\gamma}f(\omega t)dt
\tag{5.3.1}\end{equation}
Operators of a similar general type appeared in the article{\cite{Burlak1962}). J.S. Lowndes showed their connection with the Hankel transformations and separated into a indidual class (\cite {Lowndens1969}). We investigate the transformation operator for specific incoming and outgoing functions (\cite{PrudnikovBrychkovMarichev2002}, vol II, 2.12.35 (2))
\begin{equation}\label{5.3.2}
f_0(t) = t^{\nu+1}J_\nu(\omega t); \quad
f_1(x) = \beta^\mu \omega^\nu x^{\mu+\nu+1}(\beta^2+\omega^2)^{-\frac{\mu+\nu+1}{2}}J_{\mu+\nu+1}(\sqrt{\beta^2+\omega^2}x);
\tag{5.3.2}\end{equation}
and kernel
\begin{equation}\label{5.3.3}
K(x,t) = (x^2-t^2)^{\frac{\mu}{2}}J_\mu(\beta\sqrt{x^2-t^2});
\tag{5.3.3}\end{equation}
As in the previous versions, we will initially show the coefficients for which $Bf_1 (x) = 0$  and $ Af_0 (t) = 0 $, however, then more detailed attention will be paid to the eigenvalue problems. Exactly
\begin{align*}
  & a_1(x)= 1; \quad a_0(t) = 1; \\
  & b_1(x) = -\frac{1+2\mu+2\nu}{x}; \quad b_0(t) = -\frac{1+2\nu}{x}; \\
  & c_1(x) = (\beta^2+\omega^2); \quad c_0(t) = \omega^2;
\end{align*}
It is easy to verify the fulfillment of the hyperbolic equation (\ref {4.a}). For $ \mu> 0 $, the condition (\ ref {4.b}) is realized automatically, because the kernel tends to zero. If $\mu = 0 $, then $ K (0,0) = 1 $ and $ b_0 (x) = b_1 (x) $, which again confirms this equality.

As stated above, to establish the truth (\ref {4.c}), asymptotic properties of the incoming function $ f_0 (x) $ in (\ref {5.3.2}). By the definition of the Bessel function in a neighborhood of zero
\begin{equation*}
t^\nu J_\nu (\omega t) = O(t^{2\nu+1})
\end{equation*}
At the same time, the strongest singularity in (\ref {4.c}), hidden in the term with the coefficient $b_0 (\epsilon)$, has order $\varepsilon^{- 1}$. This leads to the vanishing of the entire expression (\ref {4.c}) for $\nu> 0$.

The condition at the vertex requires a cautious limit passage over each variable. The fact is that the kernel enters the root of the expression
\begin{equation*}
(\varepsilon - \delta)^2 - \delta^2 = - \delta (2\varepsilon - \delta)
\end{equation*}
The existence of an integral in the transformation operator requires the presence of an integration variable inside the segment $t \in [\epsilon, x - \delta] $. Hence, near the zero $ \delta < \varepsilon$, and this leads to the appearance of a purely imaginary factor in front of the root. The connection between the basic $ J_\mu $ and the modified Bessel functions $ I_\mu$ helps here. Exactly
\begin{equation*}
J_\mu(\imath z) = e^{\mu \frac{\pi}{2}}I_\mu(z)
\end{equation*}
Then
\begin{equation*}
K[\varepsilon, \varepsilon-\delta] = e^{\imath \mu \frac{\pi}{2}}
\sqrt{\delta (2\varepsilon - \delta)}
e^{-\imath \mu \frac{\pi}{2}}I_\mu (\sqrt{\delta (2\varepsilon - \delta)})
\end{equation*}
Exponentials are cancelled, while the remaining part for $\mu \geq 0 $ does not represent the limiting transition of difficulties. So, for (\ref {5.1.3}) all the conditions are fulfilled
(\ref {4.a} - \ref {4.d}).


\emph{Example 5.4 Generalized transmutation Vekua-Erdelyi-Lowndes operators }
\\
\begin{equation}\label{5.4.1}
(\widetilde{A}+\lambda_1)T = T(\widetilde{A}+\lambda_0)
\tag{5.4.1}\end{equation}

By definition (\cite {Sitnik2007}, \cite {RyzhkovaSitnik2016}), the Vekua-Erdelyi-Lowndes operator $\widetilde{A}$ (VEL) shifts along the spectral parameter.

We will slightly extend the formulation of the last equality, which allows us to immediately outline a rather large circle of operators of this type, leading to a similar studied class of hyperbolic Euler-Poisson-Darboux equations. The number of works on this topic is so huge that we will give only three links in which you can find a lot of sources on the presented question (\cite {ShishkinaSitnik2017}, \cite {Dettman1985}, \cite {Shishkina2018A}).
We will assume that the operator on the right repeats the not exact statement of the operator on the left side, but only its essence. In other words, the definition of transmutation will in fact be repeated
\begin{equation}\label{5.4.2}
(B+\lambda_1)T = T(A+\lambda_0)
\tag{5.4.2}\end{equation}
in which the construction of the operator B is similar to the structure of the operator A. The group properties of the transformations and the Euler-Poisson-Darboux equations begin to come to the fore here (\cite{Trimeche1988},\cite{Ibragimov2006}).

For example, consider the version of the equality (\ ref {4.a}}) of the form
\begin{align}\label{5.4.3}
  & \partial_{tt}K(x,t)-\frac{\kappa_0}{t}\partial_t K(x,t) + (\beta+\delta)K(x,t) = \notag\\
  & = \partial_{xx}K(x,t)-\frac{\kappa_1}{x}\partial_x K(x,t) + \beta K(x,t)
\tag{5.4.3}\end{align}
Naturally, to find a solution with a kernel representation
\begin{equation}\label{5.4.4}
K(x,t) = G(\Gamma); \quad \Gamma = x^2 - t^2;
\tag{5.4.4}\end{equation}
This leads to an ordinary differential equation with respect to the independent variable $\Gamma$
\begin{equation}\label{5.4.5}
4\Gamma^2 \partial_{\Gamma\Gamma}G(\Gamma) - 2(\kappa_1-\kappa_0+2)\partial_\Gamma +
(2\beta+\delta)G(\Gamma) = 0
\tag{5.4.5}\end{equation}
If for convenience we put
\begin{equation*}
\mu^2 = 2\beta+\delta; \quad  \nu = \frac{1}{2}(4-\kappa_0-\kappa_1);
\end{equation*}
then a solution of this equation will be
\begin{equation}\label{5.4.6}
K(x,t) = \Gamma^{\frac{\nu}{2}}J_\nu (\mu \sqrt{\Gamma});
\tag{5.4.6}\end{equation}

In the tables of integrals, it is easy to find an example of the transmutation (\cite{PrudnikovBrychkovMarichev2002}, vol. 2, 2.12.21 (5)), where
\begin{align}\label{5.4.7}
& f_0(t) = \cos(\omega t); \notag\\
& f_1(x) = \sqrt{\frac{\pi}{2}}x^{\nu+\frac{1}{2}}\mu^\nu (\omega^2 + \mu^2)^{-\frac{2\nu+1}{4}}J_{\nu+\frac{1}{2}}
\left(x\sqrt{\omega^2+\mu^2}\right)
\tag{5.4.7}\end{align}
It is immediately evident that the proposed kernel connects the eigenfunctions of the second derivative and the singular Bessel differential operator


\begin{equation}\label{5.4.8}
A = D^2;  \quad B = D^2 - \frac{2\nu}{x}; \quad D = \frac{\partial}{\partial x};
\tag{5.4.8}\end{equation}
with different eigenvalues.
\begin{equation}\label{5.4.9}
A f_0 = -\omega^2 f_0; \quad B f_1 = -(\omega^2 + \mu^2);
\tag{5.4.9}\end{equation}

In the special case $ \ nu = 0 $, using the well-known formulas for the Bessel function of half-integer order, a transmutation of I.N. Vekua, comparing
\begin{equation}\label{5.4.10}
f_0(t) = \cos(\omega t); \quad
f_1(x) = \frac{\sin\left(x\sqrt{\omega^2 + \mu^2}\right)}{\sqrt{\omega^2 + \mu^2}};
\tag{5.4.10}\end{equation}
The equality $f_1(x) = T f_0(t)$ can be differentiated with respect to the variable x, and then the transformation Vekua-Erdelyi-Lowndes operator $OP_ {II}$ in the 'pure' form will be constructed
\begin{equation*}
\cos\left(x\sqrt{\omega^2 + \mu^2}\right)  = \cos(\omega x)-
\mu\int\limits_0^x x\frac{J_1\left(\mu\sqrt{x^2-t^2}\right)}{\sqrt{x^2-t^2}}\cos(\omega t) dt
\end{equation*}

We note that relation
\begin{equation}\label{5.4.11}
\frac{\sin\left(x\sqrt{\omega^2 + \mu^2}\right)}{\sqrt{\omega^2 + \mu^2}} =
\int\limits_0^x J_0\left(\mu \sqrt{x^2 - t^2}\right)\cos(\omega t) dt
\tag{5.4.11}\end{equation}
is naturally a cosine-Fourier transform for the generalized function contained in the domain
$x^2 - t^2 \geq 0; \; t > 0 $, and under this category it is included in the tables (\cite {BatemanErdelyi1974}), however, in traditional summaries of integrals such as (\cite{MagnusOberhettinger1954}, Third Chapter) и (\cite{GradshteynRyshik2007}) this transformation is still included only for $\mu = 1, \, \omega = 1$ and, especially interesting for us, the variant $\omega = 0$. The importance of such a simple case is caused by the fact that the transformation operator turns the constant $ f_0 (t) = 1 $ to $ f_1 (x) = \sin (\mu t) \setminus \mu $, and the kernel is the solution of the telegraph equation
\begin{equation}\label{5.4.12}
\partial_{tt}K(x,t) - \left[\partial_{xx}K(x,t)+\mu^2 K(x,t)\right] =0;
\tag{5.4.12}\end{equation}

The group properties of VEL analogs and, often, the Euler-Poisson-Darboux equation associated with them lead (for example, \cite {Shishkina2018A}) and can lead to many more remarkable results.

\emph{Example 5.5 Variant of binding of input and output functions for equations of different order}

Consider the relation (\ cite {PrudnikovBrychkovMarichev2002}, vol I, 2.4.3 (10)), which takes $ f_0 \rightarrow f_1 $ with functions
\begin{align}\label{5.5.1} \notag
& f_0(t) = \cosh (\mu t)\\
& f_1(x) = B\left(\frac{\alpha}{2},\beta\right)
{}_1F_2\left(
\frac{\alpha}{2},\frac{1}{2}, \beta+\frac{\alpha}{2}, \left(\frac{\mu x}{2}\right)^2
\right)
\tag{5.5.1}\end{align}
where $ B (\circ, \circ) $ is a beta function. Using the guidelines on the generalized hypergeometric series (see, for example, \cite {BatemanErdelyi1974}, Special Functions, Vol. I, Ch. 4), one can check the following differential equality
\begin{equation}\label{5.5.2}
x^2\partial^{(3)}_xf_1(x) +(\alpha+2\beta)x\partial_{xx}-\mu^2 x^2 \partial_x f_1(x) -\alpha \mu^2 f_1(x) =0
\tag{5.5.2}\end{equation}
The expression on the left does not enter the domain of second-order differential operators, investigated above. However, the developed theory helps to understand this option. It suffices to note that ($ \ref {5.5.2} $) reduces to
\begin{align}\label{5.5.3}
& \partial_x\left[
x^2\partial_{xx}f_1(x)+(\alpha+2\beta-2)x\partial_x f_1(x) +(2-\alpha-2\beta-\mu^2 x^2f_1(x)
\right] - \notag\\
& -(\alpha-2)\mu^2 f_1(x) = 0
\tag{5.5.3}\end{align}
Considering the hyperbolic cosine as an input function and the possibility of invariance of the solution
\begin{equation*}
K(x,t) = G(x^2-t^2)
\end{equation*}
a hyperbolic equation with two variables of the standard form suit to it.

\begin{align}\label{5.5.4}
& \partial_{tt}K(x,t)-\frac{\alpha+2\beta}{t}\partial_t K(x,t) - \mu^2 K(x,t)
+4\frac{(\beta-1)(\beta-2)}{x^2-t^2}K(x,t) =
\notag \\
& = \partial_{xx}K(x,t) -\frac{\alpha+2\beta}{x}\partial_xK(x,t)-\mu^2 k(x,t)
\tag{5.5.4}\end{align}
It is not difficult to verify that the kernel
\begin{equation}\label{5.5.5}
K(x,t) = (x^2 - t^2)^{\beta-1}
\tag{5.5.5}\end{equation}
satisfies this equation.

\section{Concluding remarks}

The concept of determining the kernel of the transformation operator through partial differential equations originated with the appearance of transmutation. Robert Carroll, who paid much attention to the theory of transformations, repeatedly noted that between two differential expressions, must exist an integral operator connecting their solutions (\cite{Carrol1982A}). However, the implementation of this approach, has many varieties. For example, Carrol himself chose the boundless direction of conjugation of transformation operators with scattering theory(\cite{Carrol1983}). Some researchers began to consider equivalent ordinary differential equations by reproductions of the same Hilbert space with the corresponding transition kernel (Reproducing Kernel Hilbert Spaces - RKHS)(\cite{SaitohSawano2016}, ). Systems of elliptic equations whose integrals are on the basis of transmutations have been studied in (\cite{BegehrGilbert1993}).

In the work (\cite {KravchenkoTorba2018}), the search for the kernel of the integral operator was carried out by the method (\ref {4.1.1} - \ref {4.1.4}), however, not for all possible second-order equations, but only for perturbed Bessel equations.

\section*{Acknowledgments}
The author gratefully acknowledges the support of Dr. Sc. (Phys.-Math.) S. M. Sitnik.

\end{document}